# A new definition of nonlinear statistics mean and variance

W. Chen

**Present mail address** (as a JSPS Postdoctoral Research Fellow): Apt.4, West 1st floor, Himawari-so, 316-2, Wakasato-kitaichi, Nagano-city, Nagano-ken, 380-0926, JAPAN

Permanent affiliation and mail address: Dr. Wen CHEN, P. O. Box 2-19-201, Jiangshu University of Science & Technology, Zhenjiang City, Jiangsu Province 212013, P. R. China

Present e-mail: chenw@homer.shinshu-u.ac.jp

Permanent email: chenwwhy@hotmail.com

## Introduction

This note presents a new definition of nonlinear statistics mean and variance to simplify the nonlinear statistics computations. These concepts aim to provide a theoretical explanation of a novel nonlinear weighted residual methodology presented recently by the present author.

## Results

When a set of values has a sufficiently strong central tendency, it is significant to characterize the set by the concept mean. Consider the values $x_1, \ldots, x_N$, the mean of

them is defined as

$$\bar{x} = \frac{1}{N}\sum_{j=1}^{N} x_j,  \tag{1}$$

which estimates the central value of sampled point values. Variance is another useful concept to characterize the variability of this set of values, which is defined as

$$Var(x_1 \ldots x_N) = \frac{1}{N-1}\sum_{j=1}^{N}(\bar{x} - x_j)^2.  \tag{2}$$

The above equation (2) estimates the mean deviation of x from its mean value. For nonlinear variables, the same definitions of mean and variance are usually used. For example, consider a set of nonlinear values $x_1y_1, \ldots, x_Ny_N$, we have

$$\overline{xy} = \frac{1}{N}\sum_{j=1}^{N} x_j y_j,  \tag{3}$$

and

$$Var(x_1y_1 \ldots x_Ny_N) = \frac{1}{N-1}\sum_{j=1}^{N}(\overline{xy} - x_j y_j)^2.  \tag{4}$$

In this paper, we present a new definition of the mean and variance for nonlinear variable samples. Without the loss of generality, consider the above example and we have

$$\overline{\overline{xy}} = \bar{x} \bullet \bar{y} = \left(\frac{1}{N}\sum_{j=1}^{N} x_j\right)\left(\frac{1}{N}\sum_{j=1}^{N} y_j\right),  \tag{5}$$

and

$$N var(x_1y_1 \ldots x_Ny_N) = \left(\frac{1}{N-1}\sum_{j=1}^{N}(\bar{x} - x_j)^2\right)\left(\frac{1}{N-1}\sum_{j=1}^{N}(\bar{y} - y_j)^2\right).  \tag{6}$$

Also, for $sin(x_1),\ldots,sin(x_N)$, we have

$$\overline{sin(x)} = sin(\bar{x}) = sin\left(\frac{1}{N}\sum_{j=1}^{N} x_j\right),  \tag{7}$$

and

$$N\,var(sin(x_1)\ldots sin(x_N)) = sin\left(\frac{1}{N-1}\sum_{j=1}^{N}(\bar{x}-x_j)^2\right). \tag{8}$$

The same idea can be extended to the definition of nonlinear median and mode. The present definitions of the mean and variance of nonlinear data provide a statistical explanation for a novel methodology of nonlinear weighted residuals proposed in [1].

Reference:

1. Chen W., A novel methodology of weighted residual for nonlinear computations. (submitted). Also published in arXiv.org e-Print archive: http://xxx.lanl.gov/abs/math.NA/9905003 (Abstract and full text can be found there).